\documentclass[12pt]{article}
\usepackage{amsmath}
\usepackage{amssymb}

\usepackage[cp1251]{inputenc}
\usepackage[russian]{babel}
\newcommand{\il}[2]{\int\limits_{#1}^{#2}}

\newcommand{\ph}{\phantom{a}}
\newcommand{\phh}{\phantom{aaa}}

\newcommand{\sist}[2]{\left\{
\begin{array}{l}
{#1}\\
\ph\\
{#2}
\end{array}
\right.}

\begin{document}

MSC xxxxx

\vskip 15pt

\centerline{\bf Extension of Krein's special method for solving}
\centerline{\bf integral equations}
\vskip 10pt
\centerline{\bf G. A. Grigorian}

\centerline{\it Institute  of Mathematics of NAS of Armenia}
\centerline{\it E -mail: mathphys2@instmath.sci.am}
\vskip 20 pt

\noindent
Abstract. Extencion of Krein's special method for solving of integral equation to that method for solving of systems of integral equations is established. Generalizations of formulae for solution of integral equations  are obtained. The resut obtained is demostrated by examples.

\vskip 10pt

Key words: systems of integral equations, resolvent kernel, existence condition, kernel, depending on difference of arguments, symmetric kernel function.

\vskip 10pt

{\bf 1. Introduction}. Let $K(t,s)\equiv (k_{ij}(t,s))_1^n$ be a complex-valued continuous matrix function on $[a,b] \times [a,b], \ph f(t)\equiv (f_i(t))_1^n$ be a complex-valued continuous vector-function on $[a,b]$.
Consider the system of integral equations
$$
\phi(t) - \int\limits_a^b K(t,s) \phi(s) d s = f(t), \phh t \in [a,b]. \eqno (1.1)
$$
along with this equation consider the truncated equations
$$
g(t,\xi) - \il{a}{\xi}K(t,s)g(s,\xi) d s =I, \ph a \le t \le \xi\le b, \eqno (1.2)
$$
$$
g^*(t,\xi) - \il{a}{\xi}K(t,s)g^*(s,\xi) d s =I, \ph a \le t \le \xi\le b, \eqno (1.2^*)
$$
where $I$ is the $n \times n$ identity matrix.

For $n=1$ M. G. Krein proved the following theorem

\vskip 10pt

{\bf Theorem 1.1([1. p. 230, Theorem 2.30]).} {\it Let $n=1$. If for every $\xi \in(a,b]$ the truncated equation (1.2) has the unique solution $g(t,\xi) \in C(a,\xi)$ and, therefore. the truncated equation $(1,2^*)$  has the unique solution $g^*(t,\xi) \in C(a,\xi)$. Then the unique solution $\phi(t)$ of Eq. (1.1) can be found by the formula
$$
\phi(t) = \Bigl[\frac{1}{\mathcal{M}'(\xi)}\frac{d}{d\xi}
\int\limits_a^\xi g^*(s,\xi)d s\Bigr]_{x=b} g(t,b)- \int\limits_{t}^b g(t,\xi)\frac{d}{d\xi}\Bigl(\frac{1}{\mathcal{M}'(\xi)}\frac{d}{d\xi}\int\limits_a^\xi g^*(s,\xi)f(s) ds \Bigr)d\xi, \eqno (1.3)
$$
provided $\mathcal{M}'(\xi) \equiv g(\xi,\xi) g^*(\xi,\xi) \ne 0, \ph a \le \xi \le b$.
}

\vskip 10pt

In this paper we extend the Krein method for solving of scalar integral equations  (for system (1.1) with  $n=1$) to general system (1.1) (see [1,2]).

\vskip 10pt

{\bf 2.  Auxiliary relations.}

 Assume that the truncated  system
$$
\phi(t) - \int\limits_a^\xi K(t,s) \phi(s) d s = f(t), \phh (a \le t \le \xi \le b) \eqno (2.1)
$$
for every right part $f(t)$ has the unique solution. Then for every $\xi \in (a,b]$ there exists corresponding resolvent kernel, defined from the following equivalent systems
$$
\sist{\Gamma_\xi(t,s) - \int\limits_a^\xi K(t,u)\Gamma_\xi(u,s) d u = K(t,s)}
{\Gamma_\xi(t,s) - \int\limits_a^\xi \Gamma_\xi(t,u) K(u,s) d u = K(t,s), \phh (a \le t \le \xi \le b)}. \eqno (2.2)
$$
Differentiating by $\xi$ the first of these equations, multiplying at right  both sides of the second of them  by $\Gamma_\xi(\xi,s)$ we arrive to the equation (see [1, p. 226])
$$
\frac{\partial}{\partial \xi}\Gamma_\xi(t,s) = \Gamma_\xi(t,\xi)\Gamma_\xi(\xi,s). \eqno (2.3)
$$
  For any $\xi \in [a,b]$ consider the following matrix equations
$$
g(t,\xi) - \il{a}{\xi}K(t,s)g(s,\xi) d s =I, \ph a \le t \le \xi, \eqno (2.4)
$$
$$
g^*(t,\xi) - \il{a}{\xi}g^*(s,\xi)K(s,t) d s =I, \ph a \le t \le \xi, \eqno (2.4^*)
$$
where $g(t,\xi)$ and $g^*(t,\xi)$ are sought matrix-functions of dimension $n \times n$.

Assume that for every $\xi \in [a,b]$ there exists the resolvent kernel $\Gamma_\xi(t,s)$ for the system (1.1).
Using equality (2.3)
it is easy to show that the unique solutions $g(t,\xi)$ and $g^*(t,\xi)$ of equations (2.4) and $(2.4^*)$ respectively satisfy the equalities
$$
g(t,\xi) = I + \il{a}{\xi} \Gamma_\xi(t,s) d s, \eqno (2.5)
$$
$$
g^*(t,\xi) = I + \il{a}{\xi} \Gamma_\xi(s,t) d s. \eqno (2.5^*)
$$
Since $K(t,s)$ is continuous the resolvent kernel $\Gamma_\xi(t,s)$ is continuously differentiable in $\xi$. Then by (2.5) and $(2.5^*)$ the matrix functions $g(t,\xi)$ and $g^*(t,\xi)$ are continuously differentiable in $\xi$. Therefore, by differentiating (2.5) and $(2.5^*)$ and takin into account (2.3) we obtain
$$
\frac{\partial}{\partial\xi}g(t,\xi) = \Gamma_\xi(t,\xi)g(\xi,\xi), \phh (a \le t \le \xi \le b), \eqno (2.6)
$$
$$
\frac{\partial}{\partial\xi}g^*(t,\xi) = g^*(\xi,\xi)\Gamma_\xi(t,\xi), \phh (a \le t \le \xi \le b), \eqno (2.6^*)
$$

\hskip 20pt

{\bf 3. Extension of  Krein's special method}.
Consider the following matrix function
$$
M(\xi) \equiv (\xi-a) I + \il{a}{\xi}\il{a}{\xi}\Gamma_\xi(t,s)dtds.
$$
Obviously by virtue of $(2.5^*)$  and $(2.6^*)$  we have
$$
M(\xi) = \il{a}{\xi}g(t.\xi) d t = \il{a}{\xi}g^*(t,\xi)d t. \eqno (3.1)
$$
Taking into account  (2.5) or $(2.5^*)$ from here we obtain
$$
M'(\xi) = g^*(\xi,\xi) g(\xi,\xi), \phh (a \le \xi \le b). \eqno (3.2)
$$
Consequently,
$$
 M(\xi) = \il{a}{\xi}g^*(s,s) g(s,s) d s, \phh (a \le \xi \le b). \eqno (3.3)
$$

\vskip 10pt

{\bf Theorem 3.1.} {\it Let for every $\xi \in (a,b]$ equation (2.5) and, therefore equation $(2.5^*)$ have the unique solutions $g(t,\xi)$ and $g^*(t,\xi)$ respectively and $\det [g^*(\xi,\xi)g(\xi,\xi)] \ne 0, \ph a\le \xi \le b.$ Then for every continuous on $[a,b]$ vector function $f(t)$ the unique solution $\phi(t)$ of system (1.1) can be found by the formula
$$
\phi(t) = g(t,b)[M'(b)]^{-1}\Bigr[\frac{d}{d \xi}\il{a}{\xi}g^*(s,\xi) f(s) d s\Bigr]_{\xi = b} - \phantom{aaaaaaaaaaaaaaaaaaaaaaaaaaaaaaaaaaaa}
$$
$$
\phantom{aaaaaaaaaaaaaaaaaaaaaaaaa}-\il{t}{b}g(t,\xi)\frac{d}{d\xi}\Bigl\{[M'(\xi)]^{-1}\frac{d}{d\xi}\il{a}{\xi}g^*(s,\xi) f(s) d s\Bigr\} d\xi. \eqno (3.4)
$$
}

Proof. Existence of solutions $g(t,\xi)$ and $g^*(t,\xi)$ implies the existence of kernels $\Gamma_\xi(t,s)$ and $\Gamma_\xi(s,t)$ respectively  for every $\xi \in (a,b]$. Then the unique solution $\phi(t)$ of system (1.1) can be found by formulae
$$
\sist{\phi(t) = g(t) + \il{t}{b}\Gamma_\xi(t,\xi)g(\xi) d\xi,}{g(t) = f(t) + \il{a}{t}\Gamma_t(t,s)f(s) d s}  \phh (a \le t \le b). \eqno (3.5)
$$
To prove the theorem it is enough to prove equivalence (3.4) to the last system.
We set
$$
J_1(t) \equiv g(t,b)[M'(b)]^{-1}\Bigr[\frac{d}{d \xi}\il{a}{\xi}g^*(s,\xi) f(s) d s\Bigr]_{\xi = b},
$$
$$
J_2(t)\equiv -\il{t}{b}g(t,\xi)\frac{d}{d\xi}\Bigl\{[M'(\xi)]^{-1}\frac{d}{d\xi}\il{a}{\xi}g^*(s,\xi) f(s) d s\Bigr\} d\xi.
$$
Then, according to (3.4) we have
$$
\phi(t) = J_1(t) + J_2(t), \phh  a\le t \le b. \eqno (3.6)
$$
Integrating $J_2(t)$ by parts we obtain
$$
J_2(t) = - J_1(t) + g(t,t)M'(t)^{-1}\frac{d}{d t}\il{a}{t}g^*(s,t)f(s) d s + \phantom{aaaaaaaaaaaaaaaaaaaaaaaaaaaaaaaaaa}
$$
$$
\phantom{aaaaaaaaaaaaaaaaaaaa}+\il{t}{b}\Bigl(\frac{\partial}{\partial \xi}g(t,\xi)M'(\xi)^{-1}\frac{d}{d\xi}\il{a}{\xi}g^*(s,\xi) f(s) d s\Bigr) d\xi, \phh a\le t \le b.
$$
Making indicated here differentiations and taking into account (2.6), (3.2) and (3.3)   we obtain
$$
J_1(t) + J_2(t) = g(t) + \il{t}{b}\Gamma_\xi(t,\xi)g(\xi) d\xi, \phh a \le t\le b,
$$
where $g(t)\equiv f(t) + \il{a}{t}\Gamma_t(t,s)f(s) d s  \phh (a \le t \le b).$ This together with (3.5) and (3.6) implies that (3.4) is equivaleny to (3.5). The theorem is proved.

\vskip 10pt

{\bf Remark 3.1.} {System (1.1) with a continuous kernel function always reduces to a scalar integral equation, but with,  in general, a non-continuous kernel.}

\vskip 10pt

Formula (3.4) allows to find the solution of system (1.1) with "general" \hskip 2pt right part by solutions of systems (3.1) and $(3.1^*)$ with special right parts. However, the use of this formula is restricted by the condition
$$
det [g^*(\xi,\xi)g(\xi,\xi)] \ne 0, \ph a\le \xi \le b. \eqno (3.7)
$$
Below we consider a class of systems (1.1), for which the last condition is satisfied.

\vskip 10pt

{\bf 4. Systems of integral equations, with symmetric kernels, depending on difference of arguments.} Consider system (1.1) under the hypothesis that the kernel function is symmetric and depends on difference of arguments:
$K(t,s) = - H(t-s)$ and $H(t),\ph -2a \le t \le 2a,\ph$ is even.
$$
\phi(t) + \il{0}{2a}H(t-s)\phi(s) d s = f(t), \phh 0\le t \le 2a. \eqno (4.1)
$$
Then instead of equations (1.2) and $(1,2^*)$ we will have respectively
$$
g(t,\xi) + \il{a}{\xi}H(t-s)g(s,\xi) d s =I, \phh 0 \le t \le \xi, \eqno (4.2)
$$
$$
g^*(t,\xi) + \il{a}{\xi}g^*(s,\xi)H(s-t) d s =I, \phh 0 \le t \le \xi, \eqno (4.2^*)
$$
In this section we show,  that if under a general hypothesis, that $H(t)$ is integrable, then the condition (3.7) holds and formula (3.4)  for system (4.1) is valid. Assume $H(t)$ is continuous.
Since $H(t)$ is even from (4.1) we have
$$
g(\xi-t,\xi) + \il{a}{\xi}H(t-(\xi-s)g(s,\xi) d s =I, \phh 0 \le t \le \xi.
$$
Making the substitution $ s \to \xi- s$ from here we obtain
$$
g(\xi-t,\xi) + \il{a}{\xi}H(t-s)g(\xi - s,\xi) d s =I, \phh 0 \le t \le \xi.
$$
Comparing this equality with (4.2) we conclude that
$$
g(t,\xi) = g(\xi-t,\xi), \eqno (4.3)
$$
provided for every $\xi \in (0,a]$ system (4.2) has the unique solution $g(t,\xi)$. By analogy one can show that
$$
g^*(t,\xi) = g^*(\xi-t,\xi), \eqno (4.3^*)
$$
provided for every $\xi \in (0,2a]$ system $(4.2^*)$ has the unique solution $g^*(t,\xi)$
It follows from (4.2) and (4.3) that
$$
g(\xi,\xi) + \il{0}{\xi}H(s)g(s,\xi)d s = I.
$$
Differentiation of  this equality gives us
$$
\frac{d g(\xi,\xi)}{d\xi}
= -\Bigl[H(\xi) + \il{0}{\xi}H(s) \Gamma_\xi(s,\xi)ds\Bigr]g(\xi,\xi).
$$
By (2.2) from here we obtain
$$
\frac{d}{d\xi}g(\xi,\xi) = \Gamma_\xi(\xi,0)g(\xi,\xi)
$$
Since by (2.2) $g(0,0) = I$, from the last equality we obtain (the Liouville formula)
$$
\det g(\xi,\xi) = e^{\il{0}{\xi}tr \Gamma_t(0,t)dt} \neq 0, \ph \xi \in [0,2a]. \eqno (4.4)
$$
Analogously using (2.2), $(4.2^*)$ and $(4.3^*)$ one can show that
$$
\frac{d}{d\xi}g^*(\xi,\xi) = g^*(\xi,\xi)\Gamma_\xi(0,\xi).
$$
Then by the Liouwille formula we have
$$
\det g^*(\xi,\xi) = e^{\il{0}{\xi}tr \Gamma_t(0,t)^\tau dt} \neq 0, \ph \xi \in [0,2a], \eqno (4.4^*)
$$
where $\Gamma_t(0,t)^\tau$ is the transpose to   $\Gamma_t(0,t)$. This together with (4.4) implies  (3.7).
Thus we prove that if $H(t)$ is continuous and is even, then  condition (3.7) holds and formula (3.8)) is valid for system (4.1). Now we show that condition (3.7) holds for the general case, when $H(t)$ is integrable over $(a,b)$. For every $\varepsilon >0$ chose continuous on $[a,b]$ symmetric matrix fumctiom $H_\varepsilon(t)$ of dimension $n\times n$ such that $||H(t)- H_\varepsilon(t)|| \le \varepsilon$, where $||x||$ denotes a $L_1^{(n)}(a,b)$ norm of $x$. Let $g_\varepsilon(t,\xi)$ and $g_\varepsilon^*(t,\xi)$ be the unique solutionons of the integral equations
$$
g_\varepsilon(t,\xi) + \il{a}{\xi}H_\varepsilon (t-s)g_\varepsilon(s,\xi) d s =I, \ph 0 \le t \le \xi, \eqno (4.5)
$$
and
$$
g^*_\varepsilon(t,\xi) + \il{a}{\xi}g^*_\varepsilon(s,\xi)H_\varepsilon (s-t) d s =I, \ph 0 \le t \le \xi, \eqno (4.5^*)
$$
respectively and let $\Gamma_{\xi,\varepsilon}(t,s)$ be the resolvent kernel of the equation
$$
\phi(t) + \il{a}{\xi}H_\varepsilon(t-s)\phi(s) d s = f(t), \ph 0\le t \le a. \
$$
and $\Gamma_\xi(t,s)$ be the resolvent kernel of Eq. (4.1). Then, $\lim\limits_{\varepsilon \to 0}||\Gamma_{\xi,\varepsilon}(t,s) -\Gamma_\xi(t,s)|| = 0$, uniformly in $(\xi,t,s)$. On the basisis of this we chose $\varepsilon_0> 0$ so small that for every positive  $\varepsilon < \varepsilon_0$ the inequalities
$$
|tr [\Gamma_t(0,t) - \Gamma_{t,\varepsilon,t}(0,t)]| \le 1, \ph t \in [a.b]
$$
are fulfilled.
By (4.4) and $(4.4^*)$ it follows from here that
$$
\det g_\varepsilon(\xi,\xi) = e^{\il{0}{\xi}tr \Gamma_{t,\varepsilon}(0,t)dt}  =   e^{\il{0}{\xi}tr\Bigl( \Gamma_{t,\varepsilon}(0,t)dt- \Gamma_t(0,t)dt\Bigr)} e^{ \il{0}{\xi}tr \Gamma_t(0,t)dt} \ge e^{-(b-a)} e^{ \il{0}{\xi}tr \Gamma_t(0,t)dt}.
$$
$$
\det g_\varepsilon(\xi,\xi) \ge  e^{-(b-a)} e^{ \il{0}{\xi}tr \Gamma_t(t,0)^\tau dt} \ph \xi \in [0,c].
$$
unifoemly in $\xi$. Since $\lim\limits_{\varepsilon\to 0}||g(\xi,\xi) - g_\varepsilon(\xi,\xi)|| = \lim\limits_{\varepsilon\to 0}||g^*(\xi,\xi) - g^*_\varepsilon(\xi,\xi)||  = 0$ uniformly in $\xi \in[a,b]$ from the last inequalities we derive that (3.7).

Foolowing up [1] the obtained result we can formulate in the next form.

\vskip 10pt

{\bf Tjeorem 4.1} {\it Let $H(t)$ be a even matrix function from $L_1^{(n)}(a-b,b-a)$ and let for every $\xi \ph (a\le\xi\le b)$ the equation
$$
g(t,\xi) + \il{a}{\xi}H(t-s)g(s,\xi) d s = I \ph (a \le t, s ,\le \xi)
$$
has the unique continuous solution. Then $\det g^*(\xi,\xi)g(\xi,\xi) \ne 0 \ph (a\le \xi \le b)$ and  for every $f \in C^{(n)}(a,b)$  the solution $\phi$ of the equation
$$
\phi(t) + \il{a}{b}H(t - s)\phi(s) d s = f(t) \eqno (4.6)
$$
can be found by the formula
$$
\phi(t) = \Bigl[M'(\xi)^{-1}\frac{d}{d \xi}\il{a}{\xi}g(s,\xi)f(s) d s\Bigr]_{\xi= b} g(t,b) - \il{t}{b}g(t,\xi)\frac{d}{d\xi}\Bigl(M'(\xi)^{-1}\il{a}{\xi}g(s,\xi)f(s) d s\Bigr) d\xi, \ph
$$
$a \le t \le b$.
}

\vskip 10pt

{\bf Remark 4.1.} {\it Obviously  Eq. (4.6) with  integrable on $[a,b]$ kernel matrix function is always reducible to a scalar integral equation with a kermel, which in general, does not depend on  difference of arguments.}

Consider the equation
$$
q(t,\xi) + \il{-\xi}{\xi}H(t-s)q(s,\xi) ds = I, \ph (-\xi \le  t \le \xi) \eqno (4.5)
$$
We set
$$
M(\xi) \equiv \il{0}{\xi}q(s,\xi)ds \ph (0 \le \xi \le a)
$$
On the basis of the obtained above results by analogy with the proof of Theorem 8.3 from [1] (see [1, page 242]) one can prove the following theorem

{\bf Theorem 4.2.} {\it Let the even  matrix function $H(t)$ belongins to $L_1^{(n)}(-2a,2a)$, $f(t)$ belongs to  $C^{(n)}(0,2a)$ and let Eq. (4.5) has the unique solution $q(t,\xi) \ph (-\xi \le t \le \xi)$, such that $\det M'(\xi) \ne 0 \ph (-\xi \le t \le \xi)$. Then for every $f(t) \in C^{(n)}(0,2a)$ the unique solution $\phi \in C^{(n)}(0,2a)$ of Eq. (4.1) can be found by the formula
$$
\phi(t)= \frac{1}{2}M'(a)^{-1}\Bigl[\frac{d}{d t}\il{-a}{a}q(s,a)f(s) d s\Bigr]q(t,a)- \phantom{aaaaaaaaaaaaaaaaaaaaaaaaaaaaaaa}
$$
$$
-\frac{1}{2}\il{|t|}{a}q(t,\xi)\frac{d}{d\xi}\Bigl(M'(\xi)^{-1}\frac{d}{d\xi}\il{-\xi}{\xi}q(s,\xi)f(s)d s\Bigr) d\xi -\phantom{aaaaaaa}
$$
$$
\phantom{aaaaaaaaaaaaaaaaaaaaaaaaa}-\frac{1}{2}\frac{d}{d t}\il{|t|}{a}q(t,\xi)M'(\xi)^{-1}\Bigl(\il{-\xi}{\xi}q(s,\xi)d f(s)\bigr) d\xi, \ph (-a \le t \le a),
$$
where the integral $\il{-\xi}{\xi}q(s,\xi)d f(s)$ is understand in the   Stiltjes siense.
}

\vskip 10pt
{\bf Example 4.1.} {\it Assume
$$
H(t) = \begin{pmatrix} 0 & h_1(t)\\
h_2(t) & 0
\end{pmatrix}.
$$
where $h_k(t)\in L_1(-2a,2a), \ph k=1,2$. Then it is not difficult to werify that Eq. (4.5) is equivalent to the system
$$
\left\{
\begin{array}{l}
q_{11}(t,\xi) - \il{0}{\xi}h_1(t-s)ds\il{0}{\xi}h_2(s-u)q_{11}(u)d u =1,\\
q_{22}(t,\xi) - \il{0}{\xi}h_2(t-s)ds\il{0}{\xi}h_1(s-u)q_{22}(u)d u =1,\\
q_{12}(t,\xi) = - \il{0}{\xi}h_1(t-s)q_{22}(s) d s,\\
q_{21}(t,\xi) = - \il{0}{\xi}h_2(t-s)q_{11}(s) d s,\ph -\xi \le t \le \xi.
\end{array}
\right.
$$
Therefore, if $\il{-2a}{2a}|h_k(t)| d t < 1, \ph k=1,2,$ then the last system has the unique solution.Hence, Eq. (4.5) has the unique solution
$$
q(t,\xi) \equiv \begin{pmatrix} q_{11}(t,\xi) & q_{12}(t,\xi)\\
q_{21}(t,\xi) & q_{22}(t,\xi)
\end{pmatrix}, \phh 0\le t,\xi \le 2a
$$
such that $\det q(\xi,\xi) \ne 0, \ph 0 \le \xi \le 2a$.

\vskip 20pt

\centerline{\bf References}

\vskip 10pt

\noindent
1. Gohberg I C., Krein M. G., Theory of Volterra operators in Hilbert space and its \linebreak \phantom {a} applications. Moscow, ''Nauka'',  1967,  508 pages.

\noindent
2.  Krein M. G., On one method of solving linear integral equations of the first and second  \linebreak \phantom {a}  order. DAN 100 (1955), pp. 413-416.

\end{document}